\documentstyle[12pt,amssymb,amstex]{amsart}

\numberwithin{equation}{section}

\oddsidemargin
-5mm \evensidemargin -5mm \baselineskip+6pt

\def\ga{{\frak A}}
\def\gb{{\frak B}}

\def\bc{{\mathbb C}}

\def\bn{{\mathbb N}}
\def\bp{{\mathbb P}}

\def\br{{\mathbb R}}

\def\bz{{\mathbb Z}}

\def\a{\alpha}

\def\e{\epsilon}
\def\l{\lambda} 

\def\m{\mu}
\def\p{\psi}
\def\n{\nu}

\def\t{\tau}
\def\f{\varphi} 
\def\v{\phi}
\def\th{\theta}  
\def\w{\omega} 
\def\z{\zeta}

\newtheorem{thm}{Theorem}[section]
\newtheorem{lem}[thm]{Lemma}

\newtheorem{prop}[thm]{Proposition}

\newtheorem{prob}[thm]{Problem}

\def\id{\mathop{\rm id}}

\def\id{{\bf 1}\!\!{\rm I}}

\begin{document}

\small

\title[A few remarks on mixing]
{A few remarks on mixing properties of $C^*$-dynamical systems}
\author{Farrukh Mukhamedov}
\address{Farrukh Mukhamedov\\
Department of Mechanics and Mathematics\\
National University of Uzbekistan\\
Vuzgorodok, 700174, Tashkent, Uzbekistan} \email{{\tt
far75m@@yandex.ru}}
\author{Seyit Temir}
\address{Seyit Temir\\
Department of Mathematics\\
Art and Science Faculty\\
Harran University, 63200, Sanliurfa, Turkey} \email{{\tt
seyittemir67@@hotmail.com}}

\begin{abstract}
We consider strictly ergodic  and strictly weak mixing
$C^*$-dynamical systems. We prove that the system is strictly weak
mixing if and only if its tensor product is strictly ergodic,
moreover strictly weak mixing too. We also investigate some other
mixing properties of the system. \vskip 0.3cm \noindent {\it
Mathematics Subject Classification}: 46L35, 46L55, 46L51, 28D05
60J99.
\end{abstract}

\maketitle
\section{Introduction}
It is known \cite{Wa},\cite{Ru} that a notion of mixing for
dynamical systems plays an important role in quantum statistical
mechanics. A lot of papers (see, \cite{FV},
\cite{FR},\cite{L1},\cite{L2},\cite{W}) were devoted to the
investigations of mixing properties of dynamical systems. Very
recently in \cite{NSZ} certain relations between ergodicity, weak
mixing and uniformly weak mixing conditions of $C^*$-dynamical
systems  have been investigated. It is known \cite{W},\cite{KSF}
that strict ergodicity of a dynamical system is stronger than
ergodicity. Therefore, it is natural to ask,  how this notion is
related with mixing conditions. The object of this paper is to
investigate this question. Namely, we are going to consider
strictly ergodic  and strictly weak mixing $C^*$-dynamical
systems. The paper organized as follows. In section 2 we recall
some preliminaries on $C^*$-algebras and dynamical systems.
Section 3 is devoted to the characterization of strictly ergodic
$C^*$-dynamical systems. In the last section 4 we prove that the
system is strictly weak mixing if and only if its tensor product
is so. We also introduce a notion of $\phi$-ergodicity and compare
it with known mixing conditions.

\section{Preliminaries}

In this section we recall some preliminaries concerning
$C^*$-dynamical systems.

Let $\ga$ be a $C^*$-algebra with unit $\id$. An element $x\in\ga$
is called {\it self-adjoint} (resp. {\it positive}) if $x=x^*$
(resp. there is an element $y\in\ga$ such that $x=y^*y$). The set
of all self-adjoint (resp. positive) element will be denoted by
$\ga_{sa}$ (resp. $\ga_+$). By $\ga^*$ we denote the conjugate
space to $\ga$. A linear functional $\f\in\ga^*$ is called {\it
Hermitian} if $\f(x^*)=\overline{\f(x)}$ for every $x\in\ga$. A
Hermitian functional $\f$ is called {\it positive} if
$\f(x^*x)\geq 0$ for every $x\in\ga$. A positive functional $\f$
is said to be a {\it state} if $\f(\id)=1$. By $S$ (resp.
$\ga^*_h$) we denote the set of all states (resp. Hermitian
functionals) on $\ga$. Let $\ga\odot\ga$ be the algebraic tensor
product of $\ga$. By $\ga\otimes\ga$ we denote a completion of
$\ga\odot\ga$ with respect to the minimal $C^*$-tensor norm on
$\ga\odot\ga$. The set of all states on $\ga\otimes\ga$ we denote
by $S^2$.  A linear operator $T:\ga\mapsto\ga$ is called {\it
positive} if $Tx\geq 0$ whenever $x\geq 0$. A positive linear
operator $T$ is called {\it a Markov operator} if $T\id=\id$. A
pair $(\ga,T)$ consisting of a $C^*$-algebra $\ga$ and a Markov
operator $T:\ga\mapsto\ga$, is called {\it a $C^*$-dynamical
system}. In the sequel, we will call any triplet $(\ga,\f,T)$
consisting of a $C^*$-algebra $\ga$, a state $\f$ on $\ga$ and a
Markov operator $T:\ga\mapsto\ga$ with $\f\circ T=\f$, that is a
dynamical system with an invariant state, {\it a state preserving
$C^*$-dynamical system}. A state preserving $C^*$-dynamical system
is a non-commutative $C^*$-probability space $(\ga,\f)$ (see
\cite{CO}) together with a Markov operator $T$ of $\ga$ preserving
the non-commutative probability $\f$.  We say that the state
preserving $C^*$-dynamical system $(\ga,\f,T)$ is {\it ergodic}
(respectively, {\it weakly mixing, strictly weak mixing}) with
respect to $\f$ if
\begin{equation}\label{erg}
\lim_{n\to\infty}\frac{1}{n}\sum_{k=0}^{n-1}(\f(yT^k(x))-\f(y)\f(x))=0,
\ \ \textrm{for all} \ \ x,y\in\ga. \end{equation}

(respectively, \begin{equation}\label{wmix}
\lim_{n\to\infty}\frac{1}{n}\sum_{k=0}^{n-1}|\f(yT^k(x))-\f(y)\f(x)|=0,\
\ \textrm{for all} \ \ x,y\in\ga, \end{equation}

\begin{equation}\label{stmix1}
\lim_{n\to\infty}\frac{1}{n}\sum_{k=0}^{n-1}|\p(T^k(x))-\f(x)|=0,
\ \ \textrm{for all} \ \ x\in\ga, \p\in S.) \end{equation}

 The  state preserving $C^*$-dynamical system $(\ga,\f,T)$ is called {\it strictly
ergodic} with respect to $\f$ if $\f$ is the unique invariant
state under $T$.

Given a $C^*$-algebra $\ga$, by $M_n(\ga)$ we denote the set of
all $n\times n$-matrices $a=(a_{ij})$ with entries $a_{ij}$ in
$\ga$. Recall that a linear mapping $T:\ga\mapsto\ga$ is called
{\it $n$- positive} if the linear operator $T_n:M_n(\ga)\mapsto
M_n(\ga)$ given by $T_n(a_{ij})=(T(a_{ij}))$ is positive.  If $T$
is $n$-positive for all $n$ then $T$ is said to be {\it completely
positive}. It is known \cite{T} that if $T$ is a completely
positive map, then the linear operator $T\otimes
T:\ga\otimes\ga\mapsto\ga\otimes\ga$ defined by $(T\otimes
T)(x\otimes y)=Tx\otimes Ty$ is also completely positive.

\section{Strictly ergodic dynamical systems}

In this section we are going to characterize strictly ergodic
$C^*$-dynamical systems. To do it we need the following

\begin{lem}\label{inv}
Let $(\ga,\f,T)$ be  strictly ergodic. If $h\in\ga^*$ is invariant
with respect to $T$, i.e. $h(Tx)=h(x)$ for all $x\in \ga$, then
there is a number $\l\in\bc$ such that $h=\l\f.$
\end{lem}

\begin{pf} Let us first assume that $h$ is positive, then
$\p=h/h(\id)$ is a state. According to the strict ergodicity of
$(\ga,\f,T)$ we have $\p=\f$, which implies that $h=h(\id)\f$. Now
let $h$ be a Hermitian functional. Then there is a unique Jordan
decomposition \cite{T} of $h$ such that
\begin{equation}\label{her}
h=h_+-h_-, \ \ \|h\|_1=\|h_+\|_1+\|h_-\|_1,
\end{equation}
where $\|\cdot\|_1$ is the norm on $\ga^*$. The invariance of $h$
implies that
$$
h\circ T=h_+\circ T-h_-\circ T=h_+-h_-.
$$
Using $\|h_+\circ T\|_1=h_+(\id)=\|h_+\|_1$, similarly $\|h_+\circ
T\|_1=\|h_+\|_1$, from uniqueness of the decomposition we find
$h_+\circ T=h_+$ and $h_-\circ T=h_-$. Therefore, by the previous
argument one gets $h=\l\f$. If $h$ is an arbitrary functional,
then there are Hermitian functionals $h_1$,$h_2$ such that
$h=h_1+ih_2$. Again invariance of $h$ implies that $h_i\circ
T=h_i$, $i=1,2$. Consequently, we obtain that $h=\l\f$.
\end{pf}

Now we are ready to formulate a criterion for the strict
ergodicity of a dynamical system. The proof of the criterion is
similar to the proof of Theorem 2, Ch.1, sec. 8 \cite{KSF}. For
the sake of completeness we will prove it.

\begin{thm}\label{sterg}  Let $(\ga,\f,T)$ be a state preserving
$C^*$-dynamical system. The following conditions are equivalent
\begin{itemize}
\item[(i)] $(\ga,\f,T)$ is strictly ergodic ; \item[(ii)] For
every $x\in\ga$ the following equality holds
$$
\lim_{n\to\infty}\frac{1}{n}\sum_{k=0}^{n-1}T^k(x)=\f(x)\id,
$$
where convergence in norm of $\ga$; \item[(iii)] For every
$x\in\ga$ and $\p\in S$ the following equality holds
$$
\lim_{n\to\infty}\frac{1}{n}\sum_{k=0}^{n-1}\p(T^k(x))=\f(x).
$$
\end{itemize}
\end{thm}

\begin{pf} Let us consider the implication (i)$\Rightarrow$(ii).
It is clear that for every element of the form $y=T(x)-x$,
$x\in\ga$ we have
$$
\bigg\|\frac{1}{n}\sum_{k=0}^{n-1}T^k(x)\bigg\|=\bigg\|\frac{1}{n}(T^{n}(x)-x)\bigg\|\leq
\frac{2}{n}\|x\|\to 0 \ \ \ \ \textrm{as} \ \ n\to\infty.
$$
So, as $\f(y)=0$ one gets
$$
\lim_{n\to\infty}\frac{1}{n}\sum_{k=0}^{n-1}T^k(y)=\f(y)\id.
$$

It is evident that the set of elements of the form $y=T(x)-x$,\
$x\in\ga$ forms a linear subspace of $\ga$. By $\gb$ we denote the
closure of this linear subspace.  Set
$$
\gb_0=\{x\in\ga : \ \f(x)=0\}.
$$
It is clear that $\gb\subseteq\gb_0$. To show $\gb=\gb_0$ assume
that $\gb\neq\gb_0$, this means that there is an element
$x_0\in\gb_0$ such that $x_0\notin\gb$. Then according to the
Hahn-Banach theorem there is a functional $h\in\ga^*$ such that
$h\upharpoonright\gb=0$ and $h(x_0)=1$. The condition
$h\upharpoonright\gb=0$ implies that $h$ is invariant with respect
to $T$. Therefore Lemma \ref{inv} yields that $h=\l\f$, which
contradicts to $\f(x_0)=0$. Hence $\gb=\gb_0$.

Let $y\in\gb_0$. Then for an arbitrary $\e>0$ we can find
$y_\e=T(x_\e)-x_\e$ such that $\|y-y_\e\|<\e/2.$ According to the
following equality
$$
\lim_{n\to\infty}\frac{1}{n}\sum_{k=0}^{n-1}T^k(y_\e)=0
$$
there exists $n_0\in\bn$ such that
$\bigg\|\frac{1}{n}\sum\limits_{k=0}^{n-1}T^k(y_\e)\bigg\|<\e/2$
for all $n\geq n_0$. Hence, we have
\begin{eqnarray*}
\bigg\|\frac{1}{n}\sum_{k=0}^{n-1}T^k(y)\bigg\|&\leq &
\bigg\|\frac{1}{n}\sum_{k=0}^{n-1}T^k(y-y_\e)\bigg\|+\bigg\|\frac{1}{n}\sum_{k=0}^{n-1}T^k(y_\e)\bigg\|
\\
&\leq &\|y-y_\e\|+\e/2<\e \ \ \ \textrm{for all} \ \ n\geq n_0.
\end{eqnarray*}

So,
$$
\lim_{n\to\infty}\frac{1}{n}\sum_{k=0}^{n-1}T^k(y)=\f(y)\id
$$
is valid for every $y\in\gb_0$.

Now let $x\in\ga$. Put $y=x-\f(x)\id$. Obviously that $y\in\gb_0$,
and for $y$ the last equality holds, whence we get the required
relation.

The implication (ii)$\Rightarrow$(iii) is evident.  Let us prove
(iii)$\Rightarrow$(i).  Assume that $\n$ is an invariant state
with respect to $T$. According to the condition (iii) we find
$$
\lim_{n\to\infty}\frac{1}{n}\sum_{k=0}^{n-1}\n(T^k(x))=\f(x)
$$
for every $x\in\ga$. On the other hand, we have
$$
\frac{1}{n}\sum_{k=0}^{n-1}\n(T^k(x))=\n(x).
$$
Whence $\f=\n$. Thus the theorem is proved.
\end{pf}

From this Theorem we immediately infer that strict ergodicity
implies ergodicity of $C^*$-dynamical system. In the next section
we will demonstrate an  example of a dynamical system which is
ergodic but not strictly ergodic. We mention that from Theorem
\ref{sterg} one gets that  strict weak mixing trivially implies
strict ergodicity.

\section{Strictly weak mixing dynamical systems}

In this section we are going to give a criterion characterizing
strictly weak mixing $C^*$-dynamical systems.

Set
$$
\ga^*_1=\{g\in\ga^*: \ \|g\|_1\leq 1\}, \ \ \
\ga^*_{1,h}=\ga^*_1\cap\ga^*_h.
$$

Before formulating a result we recall a well known fact (see for
example \cite{Wa})

\begin{lem}\label{seq} Let $\{a_n\}$ be a bounded sequence of real
numbers. Then the following are equivalent:
\begin{itemize}
\item[(i)] $$ \lim_{n\to\infty}\frac{1}{n}\sum_{k=1}^n|a_k|=0;
$$
\item[(ii)] There exists a set $J\subset\bn$ of density zero (i.e.
$$
\lim_{n\to\infty}\frac{\textrm{cardinality}(J\cap[1,n])}{n}=0)
$$
such that $\lim\limits_{n\to\infty}a_n=0$ provided $n\notin J$;
\item[(iii)]
$$
\lim_{n\to\infty}\frac{1}{n}\sum_{k=1}^n|a_k|^2=0.
$$
\end{itemize}
\end{lem}

Now we are ready to formulate the following

\begin{thm}
\label{mix} Let $(\ga,\f,T)$ be a state preserving $C^*$-dynamical
system and $T$ be a completely positive map.  The following
conditions are equivalent:
\begin{itemize}
\item[(i)] $(\ga,\f,T)$ is strictly weak mixing;
\item[(ii)] The
state preserving  $C^*$-dynamical system
$(\ga\otimes\ga,\f\otimes\f,T\otimes T)$ is strictly weak mixing;
\item[(iii)] The
state preserving  $C^*$-dynamical system
$(\ga\otimes\ga,\f\otimes\f,T\otimes T)$ is strictly ergodic;
\item[(iv)] For every $x\in\ga$ the following equality holds
$$
\lim_{n\to\infty}\sup_{\p\in\ga^*_1}\frac{1}{n}\sum_{k=0}^{n-1}|\p(T^k(x))-\p(\id)\f(x)|=0.
$$
\item[(v)] For every $x\in\ga$ and $\p\in \ga^*$ the following equality
holds
\begin{equation}\label{stmix3}
\lim_{n\to\infty}\frac{1}{n}\sum_{k=0}^{n-1}|\p(T^k(x))-\p(\id)\f(x)|=0.
\end{equation}
\end{itemize}
\end{thm}

\begin{pf} Consider the implication (i)$\Rightarrow$(ii). Recall that complete positivity of
$T$ implies that $T\otimes T$ is so. It is clear that the state
$\f\otimes\f$ is invariant with respect to $T\otimes T$.

Let $\p,\v\in S$ be arbitrary states and $x,y\in\ker\f$. Then
according to (i) we have
$$
\lim_{n\to\infty}\frac{1}{n}\sum_{k=0}^{n-1}|\p(T^k(x))|=0, \ \ \
\ \lim_{n\to\infty}\frac{1}{n}\sum_{k=0}^{n-1}|\v(T^k(y))|=0.
$$
So according to Lemma \ref{seq} there exist two subsets
$J_1,J_2\subset\bn$ of density zero such that
\begin{equation*}\lim\limits_{n\to\infty,n\notin J_1}|\p(T^k(x))|=0, \ \ \ \ \
\lim\limits_{n\to\infty,n\notin J_2}|\v(T^k(y))|=0.
\end{equation*}
Then for the set $J=J_1\cup J_2$ we have
$$
\lim\limits_{n\to\infty,n\notin J}|\p(T^k(x))\v(T^k(y))|=0,
$$
and hence again using Lemma \ref{seq} one gets that
$$
\lim_{n\to\infty}\frac{1}{n}\sum_{k=0}^{n-1}|\p(T^k(x))\v(T^k(y))|=0.
$$
Thus,
$$\lim_{n\to\infty}\frac{1}{n}\sum_{k=0}^{n-1}|\p\otimes\v(T^k\otimes
T^k(x\otimes y))|=0.
$$

By $G$ we denote the convex hull of the set $\{\p\otimes\v: \
\p,\v\in S\}$. It is clear that the $\|\cdot\|_1$-closure of $G$
is $S^2$. Therefore given $\e>0$ and $\w\in S^2$ there is $\z\in
G$ such that $\|\w-\z\|_1<\e$. For $\z$ there is $n_0\in\bn$ such
that
$$
\frac{1}{n}\sum_{k=0}^{n-1}|\z(T^k\otimes T^k(x\otimes y))|<\e \ \
\ \ \textrm{for all} \ \ n\geq n_0.
$$
Consequently,
\begin{eqnarray}\label{eq1}
\frac{1}{n}\sum_{k=0}^{n-1}|\w(T^k\otimes T^k(x\otimes y))|&\leq &
\frac{1}{n}\sum_{k=0}^{n-1}|(\w-\z)(T^k\otimes
T^k(x\otimes y))|\nonumber \\
& & +\frac{1}{n}\sum_{k=0}^{n-1}|\z(T^k\otimes
T^k(x\otimes y))|\nonumber \\
&\leq & \|\w-\z\|_1\|x\otimes y\|+\e\nonumber \\
&<&\e(\|x\otimes y\|+1)
\end{eqnarray}
for all $n\geq n_0$.

Let $x,y\in\ga$. Denote $x^0=x-\f(x)\id$, $y^0=y-\f(y)\id$. It is
clear that $x^0,y^0\in\ker\f$. By means of \eqref{eq1}, for every
$\w\in S^2$ we have
\begin{eqnarray}\label{eq11}
\frac{1}{n}\sum_{k=0}^{n-1}|\w(T^k\otimes T^k(x^0\otimes y^0))|<\e
\ \ \ \ \textrm{for all} \ \ n\geq n_1.
\end{eqnarray}

Denote $\w_1(x)=\w(x\otimes\id), \ \w_2(x)=\w(\id\otimes x)$,
$x\in\ga$. Then according to the condition (i) there exist
$N_1,N_2\in\bn$ such that
\begin{eqnarray}\label{eq12}
\frac{1}{n}\sum_{k=0}^{n-1}|\w_1(T^k(x))-\f(x)|<\e \ \ \ \
\textrm{for
all} \ \ n\geq N_1,\nonumber \\
\frac{1}{n}\sum_{k=0}^{n-1}|\w_2(T^k(y))-\f(y)|<\e \ \ \ \
\textrm{for all} \ \ n\geq N_2.
\end{eqnarray}

Now using \eqref{eq11} and \eqref{eq12} we find
\begin{eqnarray}\label{eq2}
\frac{1}{n}\sum_{k=0}^{n-1}|\w(T^k\otimes T^k(x\otimes
y))-\f(x)\f(y)|&\leq&
|\f(y)|\left(\frac{1}{n}\sum_{k=0}^{n-1}|\w_1(T^k(x))-\f(x)|\right)\nonumber
\\
&&+|\f(x)|\left(\frac{1}{n}\sum_{k=0}^{n-1}|\w_2(T^k(y))-\f(y)|\right)\nonumber\\
&&+\frac{1}{n}\sum_{k=0}^{n-1}|\w(T^k\otimes T^k(x^0\otimes
y^0))|\nonumber \\
&<&\e(|\f(x)|+|\f(y)|+1)
\end{eqnarray}
for all $ n\geq \max\{n_1,N_1,N_2\}$.

Now let $z\in\ga\otimes\ga$. Then there exists  an element
$z_\e\in\ga\odot\ga$ such that
$$
\|z-z_\e\|<\e.
$$
It follows from \eqref{eq2} that
$$
\frac{1}{n}\sum_{k=0}^{n-1}|\w(T^k\otimes T^k(z_\e))-\f\otimes
\f(z_\e)|<\e
$$
for all $n\geq n_\e$. Therefore, we obtain
\begin{eqnarray*}
\frac{1}{n}\sum_{k=0}^{n-1}|\w(T^k\otimes T^k(z))-\f\otimes
\f(z)|&\leq & \frac{1}{n}\sum_{k=0}^{n-1}|\w(T^k\otimes
T^k(z-z_\e))|\\
& & +\frac{1}{n}\sum_{k=0}^{n-1}|\w(T^k\otimes T^k(z_\e))-
\f\otimes\f(z_\e)|\\
& & +|\f\otimes\f(z_\e-z))|\\
&\leq &\e+2\left\|z-z_\e\right\|<3\e
\end{eqnarray*}
for all $n\geq n_\e$. The last relation implies that
$(\ga\otimes\ga,\f\otimes\f,T\otimes T)$ is strictly weak mixing.

The implication (ii)$\Rightarrow$(iii) is obvious.  Let us prove
the implication (iii)$\Rightarrow$(iv). Let
$(\ga\otimes\ga,\f\otimes\f,T\otimes T)$ be strictly ergodic. Let
$x\in\ker\f, x=x^*$. Given $\e>0$, strict ergodicity of the
dynamical system (see Theorem \ref{sterg}) implies that there is
$n_{0,x}\in\bn$ such that
$$
\left\|\frac{1}{n}\sum_{k=0}^{n-1}T^k\otimes T^k(x\otimes
x)\right\|<\e \ \ \ \ \textrm{for all} \ \ n\geq n_{0,x}.
$$

Hence,
$$
\left|\frac{1}{n}\sum_{k=0}^{n-1}\p\otimes\p(T^k\otimes
T^k(x\otimes x))\right|<\e \ \ \ \ \textrm{for all} \ \ n\geq
n_{0,x},\forall\p\in\ga^*_{1,h}.
$$
As $x$ is self-adjoint we get
$$
\frac{1}{n}\sum_{k=0}^{n-1}|\p(T^k(x))|^2<\e \ \ \ \ \textrm{for
all} \ \ n\geq n_{0,x},\forall\p\in\ga^*_{1,h}.
$$
According to Lemma \ref{seq} we infer that there is
$n_{1,x}\in\bn$ such that
$$
\frac{1}{n}\sum_{k=0}^{n-1}|\p(T^k(x))|<\e \ \ \ \ \textrm{for
all} \ \ n\geq n_{1,x},\forall\p\in\ga^*_{1,h}.
$$
Consequently,
\begin{eqnarray}\label{eq3}
\sup_{\p\in\ga^*_{1,h}}\frac{1}{n}\sum_{k=0}^{n-1}|\p(T^k(x))|<\e
\ \ \ \ \textrm{for all} \ \ n\geq n_{1,x}. \end{eqnarray}

Let $x\in\ker\f$ be an arbitrary element. Then it can be
represented as $x=x_1+ix_2$, where $x_1,x_2\in\ker\f$,
$x_j^*=x_j$, $j=1,2$. It then follows from \eqref{eq3} that
\begin{eqnarray}\label{eq4}
\sup_{\p\in\ga^*_{1,h}}\frac{1}{n}\sum_{k=0}^{n-1}|\p(T^k(x))|<2\e
\ \ \ \ \textrm{for all} \ \ n\geq
n_{1,x}:=\max\{n_{1,x_1},n_{1,x_2}\}.
\end{eqnarray}

Let $\p\in\ga^*_1$. Then $\p=\p_1+i\p_2$, where
$\p_j\in\ga^*_{1,h}$,$j=1,2$. By means of \eqref{eq4} one yields
\begin{eqnarray}\label{eq5}
\sup_{\p\in\ga^*_{1}}\frac{1}{n}\sum_{k=0}^{n-1}|\p(T^k(x))|&\leq&
\sup_{\p_1\in\ga^*_{1,h}}\frac{1}{n}\sum_{k=0}^{n-1}|\p_1(T^k(x))|\nonumber\\
& & +
\sup_{\p_2\in\ga^*_{1,h}}\frac{1}{n}\sum_{k=0}^{n-1}|\p_2(T^k(x))|<4\e,
 \ \  \forall n\geq n_{1,x}.
\end{eqnarray}

Finally let $x\in\ga$. Then we have the last relation \eqref{eq5}
for the element $x^0=x-\f(x)\id$, which implies that
$$
\lim_{n\to\infty}\sup_{\p\in\ga^*_{1}}\frac{1}{n}\sum_{k=0}^{n-1}|\p(T^k(x))-\p(\id)\f(x)|=0.
$$
So the implication (iii)$\Rightarrow$(iv) is proved. The
implications (iv)$\Rightarrow$(v)$\Rightarrow$(i) are obvious.

\end{pf}

{\bf Remark.} The implication (i)$\Leftrightarrow$(v) can be
proved directly using only positivity of the operator $T$. Indeed,
it is enough to prove the implication (i)$\Rightarrow$(v). Let
$x\in\ker\f$. Assume that $\p\in \ga^*_h$ be a positive
functional. Then $\tilde{\p}(x)=\frac{1}{\p(\id)}\p(x)$ is a
state. Hence, using \eqref{stmix1} we get
\begin{equation*}
\lim_{n\to\infty}\frac{1}{n}\sum_{k=0}^{n-1}|\tilde{\p}(T^k(x))|=0
\end{equation*}
which means
\begin{equation}\label{stmix11}
\lim_{n\to\infty}\frac{1}{n}\sum_{k=0}^{n-1}|\p(T^k(x))|=0.
\end{equation}
Now let $\p\in\ga^*$ be an arbitrary functional, then it can be
represented as $\p=\sum\limits_{m=0}^3i^m\p_m$, where
$\p_m\in\ga^*_h$, ($m=0,1,2,3$) are positive functionals. By means
of \eqref{stmix11} we have
\begin{equation*}
\lim_{n\to\infty}\frac{1}{n}\sum_{k=0}^{n-1}|\p(T^k(x))|\leq
\lim_{n\to\infty}\frac{1}{n}\sum_{k=0}^{n-1}\sum_{m=0}^3|\p_m(T^k(x))|=0.
\end{equation*}
Using the same argument as in the finial part of the proof
(iii)$\Rightarrow$(iv) we obtain the required assertion.
Therefore, if we take $\p(x)=\f(yx),x\in\ga$ in (v) we easily get
\eqref{wmix}, this means that strictly weak mixing implies weak
mixing.

Using the same argument as the previous theorem one can prove the
following

\begin{thm}
\label{mix1} Let $(\ga,\f,T)$ be a state preserving
$C^*$-dynamical system and $T$ be a completely positive map.  The
following conditions are equivalent:
\begin{itemize}
\item[(i)] $(\ga,\f,T)$ is weak mixing; \item[(ii)] The state
preserving  $C^*$-dynamical system
$(\ga\otimes\ga,\f\otimes\f,T\otimes T)$ is weak mixing.
\item[(iii)] The state
preserving  $C^*$-dynamical system
$(\ga\otimes\ga,\f\otimes\f,T\otimes T)$ is ergodic.
\end{itemize}
\end{thm}

\noindent {\bf Remark.} It should be noted that  Theorem
\ref{mix1}
extends Theorem 6.3 in \cite{W} to a $C^*$-algebra setting.\\

From Lemma \ref{inv} we infer that $\f$ is a unique eigenvector
with $\l=1$ eigenvalue of multiplicity one,  for strictly ergodic
dynamical system. Now what can we say about strictly weak mixing
dynamical systems? We have the following

\begin{prop}
\label{mix2} Let $(\ga,\f,T)$ be strictly weak mixing. If there
exist a number $\a\in\bc$ with  $|\a|=1$ and $\a\neq 1$, and
$h\in\ga^*$ such that
\begin{equation}\label{inv2}
h\circ T=\a h, \end{equation} then $h=0$.
\end{prop}

\begin{pf} Assume that $h\neq 0$. Then  $h\neq \m\f$ for all
$\m\in\bc$. Then, using $|\a|=1$, one gets
\begin{eqnarray*}
\frac{1}{n}\sum_{k=0}^{n-1}|h(T^k(x))-h(\id)\f(x)|
&=&\frac{1}{n}\sum_{k=0}^{n-1}|\a^kh(x)-h(\id)\f(x)|\\
&>& \bigg||h(x)|-|h(\id)\f(x)|\bigg|>0 \ \ \ \forall n\in\bn
\end{eqnarray*}
 which contradicts to the strictly
weak mixing condition. \end{pf}

Now we are going to give a concrete example of strictly mixing
$C^*$-dynamical system.

\noindent{\bf Example 1.} Let $\ga=\ell^\infty=\{(x_n): x_n\in\bc,
\ \sup|x_n|<\infty\}$. Define an operator
$T:\ell^\infty\mapsto\ell^\infty$ by means of matrix
$(t_{ij})_{i,j\in\bn}$ such that $t_{ij}=1/2^{j}$, $i,j\geq 1$. It
is not hard to check that $\f=(1/2,1/2^2,\cdots,1/2^n,\cdots)$ is
an invariant state with respect to $T$. It is known from the
Theory of Markov Chains with countable state space (see \cite{R})
that  $\p(T^nx)$ converges to $\f$ in norm of $\ga^*$ for every
state $\p\in\ga^*$. Consequently, $T$ is strictly weak mixing.\\

The following example shows that strict ergodicity does not imply
strict weak mixing.

\noindent{\bf Example 2.} Let $S^1=\{z\in\bc: |z|=1\}$ and $\l$ be
the Lebesgue measure on $S^1$ such that $\l(S^1)=1$. Fix an
element $a=\exp(i2\pi\a)$, where $\a\in[0,1)$ is an irrational
number. Define a transformation $\t:S^1\mapsto S^1$ by $\t(z)=az$.
The measure induces a positive  linear functional
$\f_\l(f)=\int\limits_{S^1}f(z)d\l(z)$ such that $\f_\l(\id)=1$.
Consider a $C^*$-algebra $\ga=C(S^1)$, where $C(S^1)$ is the space
of all continuous functions on $S^1$.   Now by means of $\t$
define a positive linear operator $T_\t:C(S^1)\mapsto C(S^1)$ by
$(T_\t(f)(z))=f(\t(z))$ for all $f\in C(S^1)$. It is clear that
$(C(S^1),\f_\l,T_\t)$ is a state preserving $C^*$-dynamical
systems. Since $\a$ is irrational, then according Theorem 2, Ch.3
\cite{KSF} we find that the defined dynamical system is strictly
ergodic. On the other hand, it is not strictly weak mixing.
Indeed, take a linear functional $h\in C(S^1)^*$ defined by
$h(f)=\int\limits_{S^1}zf(z)d\l(z)$, $f\in C(S^1)$. Then we have
$h(T_\t(f))=a^{-2}h(f)$ for all $f\in C(S^1)$. But this
contradicts to Proposition \ref{mix2}. It should be noted that
$T_\t$ is not also weak-mixing (see \cite{Wa}, Theorem 1.27).\\

Next example shows that  strict ergodicity is stronger that
ergodicity of $C^*$-dynamical system.

\noindent {\bf Example 3.}  Consider $C^*$-algebra
$\ga=\bigotimes\limits_{\bz}M_{2} (\bc)$, where $M_{2}({\bc})$ is
the algebra of $2\times 2$ matrices over the field $\bc$ of
complex numbers. By $e^{(n)}_{ij}$, $n\in \bz, \ i,j\in\{1,2\}$ we
denote the basis matrices of the algebra $M_{2}(\bc)$ sited on
$n$th place in the tensor product
$\bigotimes\limits_{\bz}M_{2}(\bc)$. The shift automorphism
$\th:\ga\mapsto\ga$ of the algebra $\ga$ is defined by
$\th(e^{(n)}_{ij})=e^{(n+1)}_{ij}$ for every $n\in \bz$ and
$i,j\in\{1,2\}$.

Let $tr$ be the normalized trace on $M_2(\bc)$, i.e. $tr(\id)=1$.
Let ${\f_0}(\cdot)=tr(\rho(\cdot))$ be a state on  $M_2(\bc)$,
where $\rho\in M_2(\bc)$ is a positive operator  such that
$tr(\rho)=1$. Such kind of $\rho$ is called a {\it density
operator} for $\f_0$.  Now let $K:M_2(\bc)\mapsto M_2(\bc)$ be a
completely positive Markov operator such that $\f_0(x)=\f_0(Kx)$
for every $x\in M_2(\bc)$. On the algebra
$\ga_{[k,n]}=\bigotimes\limits_{[k,n]}M_2(\bc)$ define the
following linear functional
$$
\f_{[k,n]}(a_k\otimes a_{k+1}\otimes...\otimes
a_n)=\f_0(a_kK(a_{k+1}(\cdots K(a_n)\cdots))).
$$
The defined functional $\f_{[k,n]}$ is a state (see
\cite{AcF},\cite{AcL}). If a compatibility condition holds
$$
\f_{[k,n]}\upharpoonright\ga_{[k-1,n-1]}=\f_{[k-1,n-1]}
$$
for the states $\{\f_{[k,n]}\}$, then there is a state $\f_K$ on
$\ga$ such that $\f_K\upharpoonright\ga_{[k,n]}=\f_{[k,n]}$ (see
\cite{Ac}), and $\f$ is called a {\it Markov state}.  We note that
a more general definition of Markov state was given in
\cite{Ac},\cite{AcF}.

It is easy to see that the Markov state is invariant with respect
to $\th$. Define two Markov operators $K_i:M_2(\bc)\mapsto
M_2(\bc)$, $i=1,2$ by
$$ K_1\left(
\begin{array}{cc}
a & b \\
c & d \\
\end{array}
\right)= \left(
\begin{array}{cc}
p_{11}a+p_{12}d & 0\\
0 & p_{21}a+p_{22}d \\
\end{array}
\right),  \ \ \ \ K_2\left(
\begin{array}{cc}
a & b \\
c & d \\
\end{array}
\right)= \left(
\begin{array}{cc}
q_{1}a+q_{2}d & 0\\
0 & q_{1}a+q_{2}d \\
\end{array}
\right).
$$

Here $\bp=(p_{ij})$ is a stochastic matrix, such that $p_{ij}>0$
for all $i,j$, and $q_1+q_2=1$, $q_1,q_2>0$.

Now consider two states $\f_{0,1}$ and $\f_{0,2}$ defined on
$M_2(\bc)$, whose density operators are given by
$$
\rho_1=
\left(\begin{array}{cc}
p_1 & 0 \\
0 & p_2 \\
\end{array}
\right ), \ \ \ \rho_2= \left(\begin{array}{cc}
q_1 & 0 \\
0 & q_2 \\
\end{array}
\right )
$$
where $\pi=(p_1,p_2)$ is a vector such that  $p_1+p_2=1$, $p_1\geq
0, \ p_{2}\geq 0$ and  $\pi \bp=\pi$.

Note that for these operators and states the compatibility
condition is satisfied, therefore there are two  associated Markov
states $\f_{K_1}$ and $\f_{K_2}$.

Then $(\ga,\f_{K_1},\th)$ and $(\ga,\f_{K_2},\th)$ are weak
mixing, and hance ergodic, state preserving $C^*$-dynamical
systems (see \cite{GN}, Th-ms 4.1 and 4.5). On the other hand,
they are not strictly ergodic,
because there exist two invariant states with respect to $\th$.\\

\noindent{\bf Remark.} From Examples 2 and 3 we conclude that weak
mixing and strict ergodicity are not comparable. Therefore we may
formulate the following

\begin{prob} Let $(\ga,\f,T)$ be a state preserving
$C^*$-dynamical system. Are the following conditions equivalent?
\begin{itemize}
\item[(i)] $(\ga,\f,T)$ is weak mixing and strictly ergodic;
\item[(ii)] $(\ga,\f,T)$ is strictly weak mixing.\\
\end{itemize}
\end{prob}

Recall a state preserving dynamical system $(\ga,\f,T)$ is called
{\it exact} (see \cite{L1}) if, for each $\p\in\ga^*$,
\begin{equation*}
\lim_{n\to\infty}\|\p\circ T^n-\p(\id)\f\|_1=0
\end{equation*}
is valid, where $\|\cdot\|_1$ is the norm in $\ga^*$. It is not
hard to see that the exactness implies  strict weak mixing. In
\cite{L1} $\L$uczak proved that  exact and weak mixing conditions,
for dynamical semi-groups on von Neumann algebras, are equivalent
if and only if the von Neumann algebra is strongly $\br_+$-finite.
Regarding this result we can formulate the following

\begin{prob} Let $(\ga,\f,T)$ be a state preserving
$C^*$-dynamical system. When are the following conditions
equivalent?
\begin{itemize}
\item[(i)] $(\ga,\f,T)$ is exact;
\item[(ii)] $(\ga,\f,T)$ is strictly weak mixing.\\
\end{itemize}
\end{prob}

Now by $S_0$ denote the set of all continuous functionals
$f:\ga_{+}\mapsto\br_+$ such that
$$
f(\l x)=\l f(x) \ \ \ \textrm{for all} \ \ \l\in\br_+,\
x\in\ga_{+},
$$
$$
f(\id)=1.
$$

Now we introduce  a notion of $\phi$-ergodicity. Namely,  a state
preserving $C^*$-dynamical system $(\ga,\f,T)$ is called {\it
$\v$-ergodic} if the equality
\begin{equation}\label{f-erg}
f(T(x))=f(x) \ \ \ \textrm{for all} \ \ x\in\ga_{+},
\end{equation}
where $f\in S_0$, implies that $f(x)=\f(x)$ for all $x\in\ga_+$.

\begin{thm}
\label{mix3} Let $(\ga,\f,T)$ be a state preserving
$C^*$-dynamical system.  Then for the conditions:
\begin{itemize}
\item[(i)] For every $x\in\ga$ the following equality holds
$$
\lim_{n\to\infty}\frac{1}{n}\sum_{k=0}^{n-1}\|T^k(x)-\f(x)\id\|=0;
$$
\item[(ii)] For every $x\in\ga$ the following equality holds
$$
\lim_{n\to\infty}\|T^n(x)-\f(x)\id\|=0;
$$
\item[(iii)] $(\ga,\f,T)$ is $\v$-ergodic; \item[(iv)] For every
$\p\in S$ and $x\in\ga$ the following equality holds
$$
\lim_{n\to\infty}\p(T^n(x))=\f(x)
$$
\end{itemize}
the following implications hold: (i)$\Leftrightarrow
$(ii)$\Rightarrow$(iii)$\Rightarrow$(iv)
\end{thm}

\begin{pf} The (i)$\Leftarrow$(ii) implication is obvious.
Consider the implication (i)$\Rightarrow$(ii). Assume that
$x\in\ker\f$, then we have
\begin{equation}\label{f-erg1}
\lim_{n\to\infty}\frac{1}{n}\sum_{k=0}^{n-1}\|T^k(x)\|=0.
\end{equation}
On the other hand, one gets
$$
\|T^{n+1}(x)\|\leq\|T^n(x)\|
$$
this means that the sequence $\{\|T^nx\|\}$ is non-increasing.
Hence, we have $\lim\limits_{n\to\infty}\|T^nx\|=\a$. It follows
from \eqref{f-erg1} that $\a=0$.  Let $x\in\ga$, then setting
$x^0=x-\f(x)\id$ we find
$$
\lim_{n\to\infty}\|T^n(x^0)\|=0
$$
which implies (ii).

(ii)$\Rightarrow$(iii). Assume that \eqref{f-erg} is valid for
some $f\in S_0$. According to the condition  (ii) we have
$$
T^n(x)\to\f(x)\id \ \ \ \ \textrm{as} \ \ \ n\to\infty
$$
for $x\in\ga_{+}$, here the convergence in norm of $\ga$. By means
of continuity of $f$ one gets
$$
f(T^n(x))\to f(\f(x)\id)=\f(x) \ \ \ \textrm{as} \ \ \ n\to\infty.
$$
On the other hand \eqref{f-erg} implies that $f(x)=\f(x)$,
$\forall x\in\ga_+$. So $(\ga,\f,T)$ is $\v$-ergodic.

(iii)$\Rightarrow$(iv). Let $\p\in S$. Define functionals $\hat
f:\ga_{+}\mapsto\br_+$, $\check f:\ga_{+}\mapsto\br_+$ by
$$
\hat f(x)=\limsup_{n\to\infty}\p(T^n(x)), \ \ x\in\ga_{+},
$$
$$
\check f(x)=\liminf_{n\to\infty}\p(T^n(x)), \ \ x\in\ga_{+}.
$$
It is clear that $\hat f,\check f\in S_0$. We have
$$
\hat f(Tx)=\limsup_{n\to\infty}\p(T^{n+1}(x))=\hat f(x).
$$
Similarly $\check f(Tx)=\check f(x)$. Hence, $\v$-ergodicity of
$(\ga,\f,T)$ implies that
$$
\hat f(x)=\f(x), \ \ \check f(x)=\f(x) \ \ \ \ \forall x\in\ga_+.
$$
Consequently, we infer the existence of the following limit
$$
\lim_{n\to\infty}\p(T^n(x))=\f(x), \ \ \ x\in\ga_{+}.
$$
Every $x\in\ga$ can be written as $x=\sum\limits_{m=0}^3i^mx_m$,
$x_m\in\ga_{+}$, ($m=0,1,2,3$), therefore by means of the last
equality we get
$$
\lim_{n\to\infty}\p(T^n(x))=\f(x), \ \ \ x\in\ga.
$$

This completes the proof.
\end{pf}

From this theorem we have

\begin{prob} Is the implication (iv)$\Rightarrow$(iii)  true?
\end{prob}

It is clear that the exactness of a dynamical system implies the
condition (iv). Therefore, it is natural to formulate the
following

\begin{prob} Let $(\ga,\f,T)$ be a state preserving
$C^*$-dynamical system. How are the following conditions related
with each other?
\begin{itemize}
\item[(i)] $(\ga,\f,T)$ is $\phi$-ergodic;
\item[(ii)] $(\ga,\f,T)$ is exact.
\end{itemize}
\end{prob}

\noindent {\bf Remark.} If $C^*$-algebra $\ga$ is finite
dimensional then all conditions in Theorem \ref{mix3} are
equivalent.\\

\section*{acknowledgements}

The author (F.M.) acknowledges the TUBITAK-NATO PC- B programme
for providing financial support, and Harran University for all
facilities they provided and kind hospitality. Authors also thank
with gratitude to a referee for useful observations and
suggestions.


\begin{thebibliography}{9999}

\bibitem{Ac} Accardi, L.,  On noncommutative Markov property,
Funct. Anal. Appl. {\bf 9} (1975), 1--7.

\bibitem{AcF} Accardi, L. and Frigerio, R., Markovian cocycles, Proc.
R. Ir. Acad. {\bf 83A} (1983), 251--263.

\bibitem{AcL} Accardi, L. and  Liebscher, V., Markovian KMS-states for
one-dimensional spin chains, Infin. Dimens. Anal. Quantum Probab.
Relat. Top. {\bf 2} (1999), 645--661.

\bibitem{CO} Cuculescu, I. and Oprea, A.G., Noncommutative probability,  Kluwer AP,
Dordrecht, 1994. .

\bibitem{FV} Frigerio, A. and Verri, M., Long-time asymptotic properties of dynamical semi-groups
on $W^*$-algebras. Math. Z. {\bf 180}(1982), 275-286.

\bibitem{FR} Fagnola, F. and Rebolledo R., Transience and recurrence of quantum Markov semi-groups.
Probab. Theory Relat. Fields {\bf 126}(2003), 289-306.

\bibitem{GN} Golodets, V.Ya. and Neshveyev, S.V., Non-Bernoullian
quantum K-systems. Commun. Math. Phys. {\bf 195}(1998), 213-232.

\bibitem{KSF} Kornfeld, I.P., Sinai, Ya. G. and Fomin, S.V., {\it Ergodic Theory},  Springer,
Berlin--Heidelberg--New York, 1982.

\bibitem{L1} $\L$uczak, A.,  Quantum dynamical semi-groups in
strongly finite von Neumann algebras, Acta Math. Hungar. {\bf
92}(2001), 11-17.

\bibitem{L2} $\L$uczak, A., Mixing and asymptotic properties of Markov
semi-groups on von Neumann algebras, Math. Z. {\bf 235}(2000),
615-626.

\bibitem{NSZ} Nicolescu, C., Str\"oh, A. and Zsid\'o, L., Noncommutative
extension of classical and multiple recurrence theorems,
J.Operator Theory, {\bf 50}(2003), 3-52.

\bibitem{R} Revuz D. {\it Markov chains},
North--Holland, Amsterdam 1984.

\bibitem{Ru} Ruelle, D., {\it Statistical mechanics},
Benjamin, Amsterdam--New York, 1969.

\bibitem{T} Takesaki, M., {\it Theory of Operator algebras, I}, Springer,
Berlin--Heidelberg--New York, 1979.

\bibitem{Wa} Walters, P., {\it An introduction to ergodic theory},  Springer,
Berlin--Heidelberg--New York, 1982.

\bibitem{W} Watanabe, S., Asymptotic behaviour and eigenvalues of dynamical semigroups
on operator algebras, Jour. Math. Anal. Appl. {\bf 86}(1982),
411-424.

\end{thebibliography}
\end{document}